\documentclass{article}

\usepackage{helvet}
\usepackage{courier}
\usepackage{graphicx}
\usepackage{amsmath}
\usepackage{float}
\usepackage{graphicx}
\usepackage{wrapfig}

\usepackage{breqn}
\usepackage{amsfonts}
\usepackage{amssymb}

\newtheorem{prob1}{Open Problem}

\newtheorem{definition}{Definition}

\newtheorem{theorem}{Theorem}
\begin{document}

\title{Exploring Tetris as a Transformation Semigroup}

\author{Peter C. Jentsch, University of Waterloo\\
        Chrystopher~L.~Nehaniv,\\ Waterloo Algebraic Intelligence \& Computation Laboratory,\\  University of Waterloo
}

\date{Received: date / Accepted: date}

\maketitle

\abstract{
    Tetris is a popular puzzle video game, invented in 1984. We formulate two versions of the game as a transformation semigroup and use this formulation to view the game through the lens of Krohn-Rhodes theory. In a variation of the game upon which it restarts if the player loses, we find permutation group structures, including the symmetric group $S_5$ which contains a non-abelian simple group as a subgroup. This implies, at least in a simple case, that iterated Tetris is finitarily computationally universal.
}
\section{Introduction}
\label{intro}

\begin{figure}
    \centering
    \includegraphics[width = 200px]{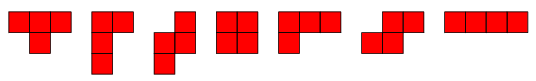}
    \caption{Pieces available in standard Tetris, the one-sided tetrominoes \cite{weisstein2003}.}
    \label{Tetris_pieces}
\end{figure}

Tetris is an arcade puzzle game created by Alexey Pajitnov in 1984, that has since become a worldwide cultural phenomenon \cite{hoad_2014}. It is the best selling paid-downloaded mobile game of all time, with over 100 million copies sold for cellphones \cite{GWR_downloaded}. It is also the most ported video game ever, according to the Guinness Book of World Records, with an estimated 65 platforms \cite{GWR_downloaded}.
Tetris is fundamentally a polyomino stacking game. The playing field consists of a 10 $\times$ 20 grid, and the player is given a sequence of tetrominoes (Figure \ref{Tetris_pieces}), which are sets of four connected grid cells, to drop from the top of the playing field. The player can translate or rotate the shapes as they fall. If a row is filled, the row disappears, and all the blocks above that row are moved down by one row. If the blocks are stacked outside the grid, then the game is over. The object of the game is to survive as long as possible. While generally, the player is only allowed to see one or two pieces ahead, most authors consider the version where the player has access to the full sequence of pieces ahead of time. This variation is also called \textit{offline} Tetris \cite{demaine2003Tetris}, and unless otherwise stated will be the one discussed here.

There has been considerable research into the mathematics behind Tetris, for an arcade game. Demaine et al.\ show that the complexity of solving many aspects of the game, such as maximizing the number of rows cleared, or the number of moves before the game ends, are NP-complete. Furthermore, they show that finding algorithms which approximate solutions to these is quite difficult \cite{demaine2003Tetris}. Other authors have characterized optimal strategies for small subsets of the pieces \cite{brzustowski1992can}, and characterized sequences of pieces that always cause a loss \cite{burgiel1997lose}.  Hoogeboom and Kosters~\cite{hoogeboom2005theory} show that nearly any reasonable configuration of blocks is possible to construct under the Tetris rules with a suitable sequence of tetrominoes. It is also possible to represent Tetris, and other tiling games, as regular grammars, which has allowed for some enumeration of possible Tetris games \cite{baccherini2008combinatorial}.

We formulate the game of Tetris as a transformation semigroup, where the elements of the semigroup are transformations on the set of possible game states. Krohn-Rhodes theory~\cite{krohn1965algebraic,maler2010krohn} and the related holonomy decomposition (Theorem~\ref{holonomy}) \cite{eilenberg1976automata,egri2017ideas} provide a way to decompose transformation semigroups into wreath products of finite simple groups and the flip-flop monoid (see Appendix for concepts and theorems employed here related to holonomy). Our analysis is primarily computational, and we use a package for the computer algebra system GAP called ``SgpDec" \cite{egri2014sgpdec}.

\section{Tetris as a transformation semigroup}

Let $P$ be a set of pieces, where a ``piece" is a set of connected cells such as a tetromino. Let $S$ be a semigroup generated by basic events $\sigma = (p,\xi) \in S$ consisting of a set of connected cells $p$, and a position $1\leq \xi \leq n$ (although the precise limits on the position $\xi$ depend on the width of $p$).

A configuration, or state, $x$ is an element of the set of $n \times k$ board of cells, where some of the cells are filled by other pieces. An element $\sigma \in S$ acts on a configuration $x$ by ``dropping" the piece $p$ with the \emph{leftmost} block in the column $\xi$, and then applying the row removal rules in the familiar way. That is, if there is a full row of width $n$, it is removed, and the blocks above the row are dropped down by one. (If more full rows arise, removal and dropping is iterated.) We denote the empty game state before any pieces have dropped by $e$, and denote by $E$ the ``game over'' state. If the total number of cells in the stack exceeds $k$, then $x \cdot \sigma = E$. Furthermore, $E \cdot \sigma = E$ for all $\sigma \in S$. For $\sigma_1, \sigma_2 \in S$, define their product $\sigma_1\sigma_2$ as the transformation resulting from applying $\sigma_1$ then $\sigma_2$ in the above way. 

We consider a state $x$ in the set of possible permutations of a $n \times k$ board of cells to be ``reachable'' in $S$ if it can be constructed by playing the game from an empty board with some sequence of pieces. More precisely, a state $x$ is ``reachable'' if there exists a word $$\sigma_1\sigma_2... \sigma_i \in S$$ such that $$x = e \cdot \sigma_1\sigma_2... \sigma_i.$$ The semigroup is precisely the set of transformations $S$ given by concatenating pieces and possible positions for those pieces on the board $(p,\xi)$, along with the set of game states reachable from the empty board using those $(p,\xi) \in S$.

\begin{definition}
$(X,S)$ is a finite transformation semigroup, which we will call the \textbf{Tetris semigroup} of $P$ on the board with dimensions $n \times k$.
\end{definition}

\section{Analysis}

In order to obtain a full description of $(X,S)$ in terms of transformations, we implemented the rules of Tetris in Python. Here, we will consider the complexity of a few variants of the game. Computation can be done on the $3 \times 3$ and $3 \times 4$ size gameboard with tri-ominoes as described below, but any larger is currently out of reach of the computational capabilities of the GAP algebra system and SgpDec.

\subsection{Tri-tris}

\begin{figure}[H]
    \centering
    \includegraphics[width = 200px]{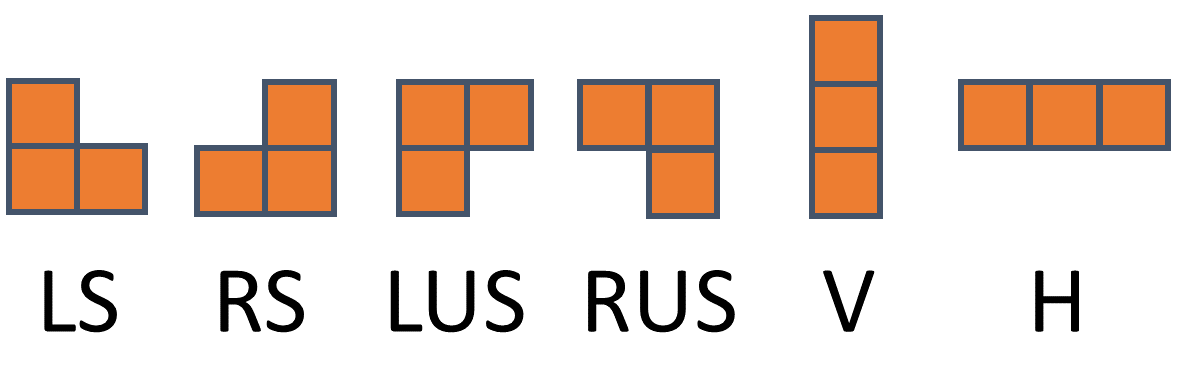}
    \caption{Triominoes and corresponding labels used in the implementation.}
    \label{triominoes}
\end{figure}

Standard Tetris has an extremely large state-space. Germundsson \cite{germundsson1991Tetris} estimates that it is on the order of $2^{200}$, this estimate is corroborated by the later constructibility result of \cite{hoogeboom2005theory}. Therefore we will consider a variant of Tetris on an $n \times k$ board, using \textit{triominoes} (Figure \ref{triominoes}) rather than tetrominoes.

Let $P = \{\text{LS},\text{RS},\text{LUS},\text{RUS},\text{V},\text{H}\}$, then the corresponding game (and semigroup) we will call \emph{Tri-tris} accordingly.

We will mostly consider the setting where $n = 3$. In this case, the horizontal line triomino $H$ is the identity, except when there are cells in the top line, in which case it maps to $E$. We will exclude $H$ as it doesn't add significantly to the system.

\subsection{Aperiodicity}
A transformation semigroup is called {\em aperiodic} if all of its subgroups are trivial. We have found Tetris to be aperiodic for $n \leq 3$ and $k \leq 5$. We pose this as an open problem for larger board sizes. 

\begin{prob1}[Is Tetris aperiodic?]
For any $\sigma \in S$ does there exist a $k>0$ such that, for all $x\in X$, $x \cdot \sigma^{k+1} = x \cdot \sigma^k$. If so, the Tetris semigroup $(X,S)$ is aperiodic, and the corresponding KR decomposition contains no nontrivial permutation groups.
\end{prob1}

The aperiodic complexity of a transformation semigroup is the least number of identity-reset components (i.e., direct products of flip-flops) that must be wreathed together to emulate it. 

\begin{table}[H]
    \centering
    \begin{tabular}{l  l  c l}
    
\hline\noalign{\smallskip}
        Board Dimensions & $|X|$ & $|S|$ & $h_s(X)$\\

\noalign{\smallskip}\hline\noalign{\smallskip}
        $3 \times 3 $ & 35 & 2,056 & 13 \\
        $3 \times 4 $ & 135 & 259,726 & 32 \\
        $3 \times 5$ & 709 & - & - \\
        
\hline\noalign{\smallskip}
    \end{tabular}
    \caption{Bounding the aperiodic complexity of Tri-tris on boards of different size by the length of the longest proper subduction chain, $h_s(X)$. All boards with width $n = 3$ have 11 generators. When $n$ is increased to $4$ the semigroup is too large for GAP to handle.}
    \label{aperiodic_standard}
    
\end{table}

If Tetris is always aperiodic, this means that there are no internal symmetries for the holonomy decomposition to expose. The complexity increases extremely quickly as the board size increases. In the next section, we will consider a rule modification that introduces these symmetries.

\section{Periodic Tri-tris}

The most straightforward rule modification to Tetris that gives the system reversibility (groups in the holonomy decomposition) is replacing the end state $E$ with the empty board $e$. In this version of the game, any move that would previously have caused a loss, now returns the game to the empty board.

\begin{table}[H]
    \centering
    \begin{tabular}{llcl}
    
\hline\noalign{\smallskip}
        Board Dimensions & $|X|$ & $|S|$ & holonomy groups present\\
        
\noalign{\smallskip}\hline\noalign{\smallskip}
        $3 \times 3 $ & 34 & 118,637 & $(4,C_2 \times C_2), (3,S_3), (2,C_2)$ \\
        $3 \times 4 $ & 135 & - & - \\
        $3 \times 4$, $P = \{\text{RS},\text{LUS},\text{RUS}, \text{V}\}$ &  116 & - &  $(4,C_2), (5,S_5), (4,S_4), (3,S_3), (2,C_2)$\\
\hline\noalign{\smallskip}
    \end{tabular}
    \caption{The complexity of periodic Tri-tris on boards of different sizes, showing the groups present in the holonomy decomposition given by SgpDec. Replacing the losing state $E$ with the empty board state $e$, the semigroups become much larger. A computation for a $3 \times 4$ board with a reduced generator set is included as well.}
    \label{aperiodic_standard}
\end{table}

\subsection{Periodic Tri-tris: $3 \times 3$ Case}

\begin{figure}
    \centering
    \includegraphics[width =\textwidth]{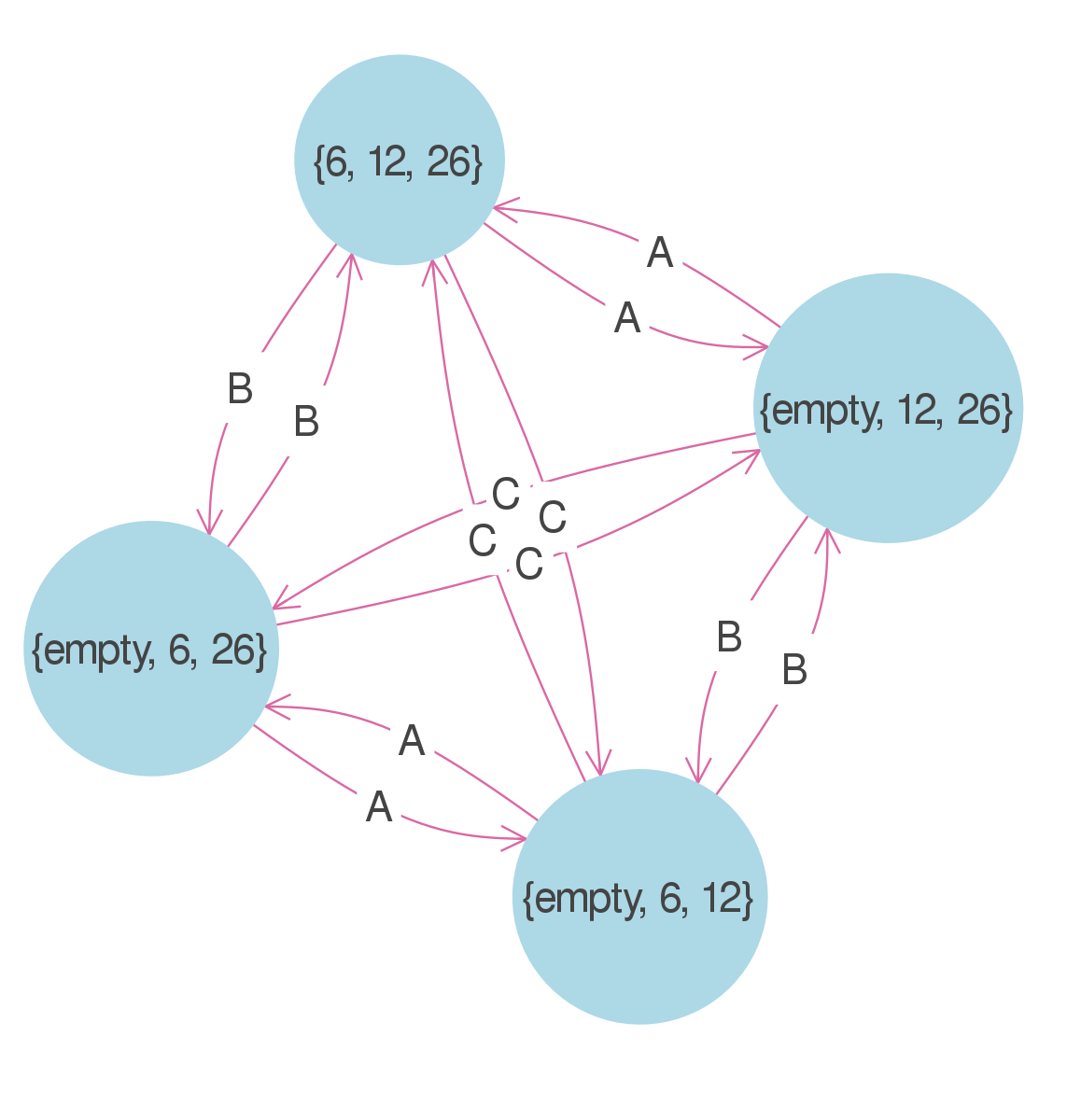}
    \caption{Holonomy transformation group $(4,C_2 \times C_2)$ where $A = V_0 LS_1 V_2 V_0 RS_1 V_1 V_2 V_0 RS_1 V_1$, $B = V_0 LS_1 V_2 V_1 V_0 V_2 LS_0 V_1 V_0 V_2 LS_0 V_1$, $C = V_0 V_2 LS_0 V_1$}
    \label{c2xc2}
\end{figure}

 Non-trivial holonomy permutation groups appear in this case.  Consider the empty state together with the three other states: $\{\mbox{empty}, 6, 12, 26\}$, whose non-empty states are visualized pictorially in Figure~\ref{states}.  The  Figure~\ref{c2xc2} shows how the members of the holonomy group $C_2 \times C_2$ act on tiles of this set. We can find pictorial representations of some of the states and the transformations to more easily visualize what this figure is describing within the games (Figure \ref{states}). Generally, the way that these groups permute the tile is that the words ``reset" the state by exceeding the length of the board, and then construct the new state. The non-abelian group $S_3$ also appears in the holonomy decomposition.

\begin{figure}
    \centering
   \includegraphics[width=0.3\textwidth]{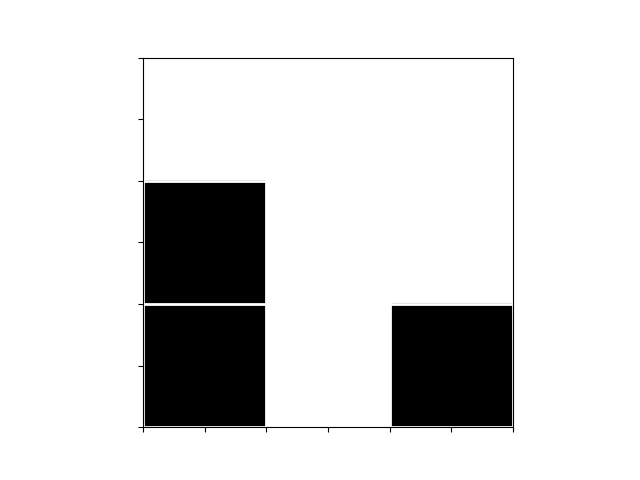}\label{state6}
    \caption{Visualizations of state 6 in Figure \ref{c2xc2}.}
    \label{states}
\end{figure}

\begin{figure}
    \centering
   \includegraphics[width=0.3\textwidth]{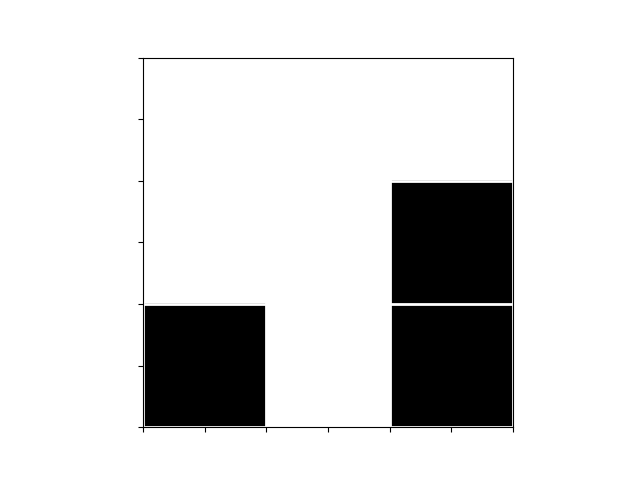}\label{state12}
    \caption{Visualizations of state 12 in Figure \ref{c2xc2}.}
    \label{states}
\end{figure}

\begin{figure}
    \centering
   \includegraphics[width=0.3\textwidth]{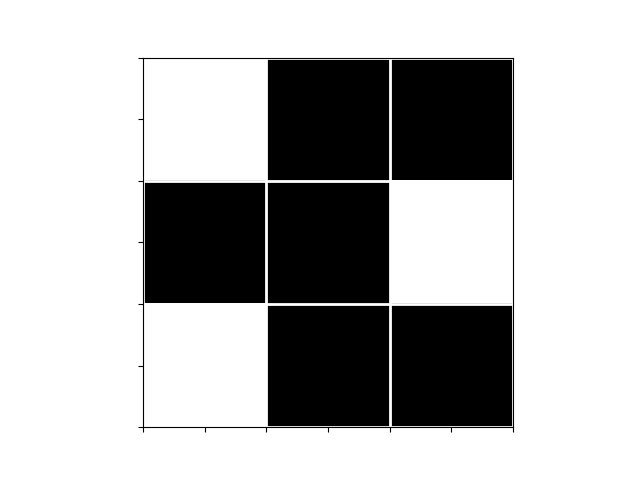}\label{state26}
    \caption{Visualizations of state 26 in Figure \ref{c2xc2}.}
    \label{states}
\end{figure}

\subsection{Periodic Tri-Tris: $3 \times 4$ Case with Reduced Generating Set} 
\begin{figure}[h]
    \centering
    \includegraphics[width = \textwidth]{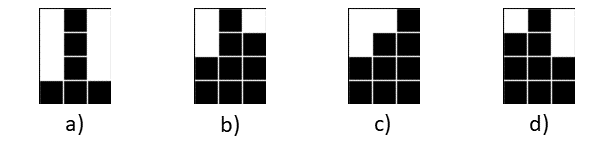}
    \caption{Visualizations of the states permuted by $S_5$, in addition to the empty state, in the Tetris semigroup with $n \times k = 3 \times 4$, $P = \{RS,LUS,RUS,V\}$. a) State 4, b) State 11, c) State 13, d) State 16.}
    \label{states_s5}
\end{figure}

If we increase the board size to $3 \times 4$, and let $P = \{RS,LUS,RUS,V\}$, we see that the holonomy decomposition contains the full symmetric group $S_5$, acting on the set $Z = \{\text{empty}, 4, 11, 13, 16\}$, illustrated in Figure \ref{states_s5}. Although this group is too large to visualize in the same fashion as $C_2 \times C_2$, we can describe some generators of the permutator group. For instance, the word in equation \ref{gen1} is a 5-cycle on $Z$, and the word in equation \ref{gen2} is a 2-cycle on $(11,16)$ and a 3-cycle on $(\text{empty},4,13)$. \\

\begin{dmath}
\label{gen1}
    V_1LUS_0V_2V_0RS_1RS_1LUS_0RS_1V_0V_1V_2RS_0RS_0V_2V_0V_1V_0V_2V_1V_2LUS_0RS_1\\V_0V_1V_2RS_0RS_0V_2V_0V_1V_0V_2V_1V_2LUS_0RS_1V_0V_1V_2RS_0RS_0V_2V_0V_1V_0V_2V_1V_2\\LUS_0RS_1V_0V_1V_2RS_0
    RS_0V_2V_0V_1V_0V_2V_1V_2\\LUS_0RS_1V_0V_1V_2RS_0RS_0V_2V_0V_1V_0V_2V_1V_2LUS_0RS_1V
_0V_1V_2RS_0RS_0V_2V_0\\ = (\text{empty},4,11,13,16)
\end{dmath}
\begin{dmath}
\label{gen2}
    V_1LUS_0V_2V_0RS_1RS_1V_0V_1RS_0V_2V_0V_2LUS_1V_0V_1V_2RS_0RS_0V_2V_0V_1RUS_1\\
    V_0V_2LUS_1V_0V_1V_2RS_0RS_0V_2V_0V_1RUS_1V_0V_2LUS_1V_0V_1V_2RS_0RS_0V_2V_0V_1\\
    RUS_1V_0V_2LUS_1V_0V_1V_2RS_0RS_0V_2V_0  = (\text{empty},4,13)(11,16)
\end{dmath}
These two permutations generate the group $S_5$.

\section{SNAGs and Computation in Tetris}
The appearance of the symmetric group $S_5$ in the transformation semigroup of 3$\times$4 periodic Tri-tris shows that the smallest simple nonabelian group $A_5$ (the alternating group of 60 even permutations on 5 elements) can be emulated by this Tetris semigroup. Ensemble techniques of Nehaniv et al.\ \cite{nehaniv2015symmetry}
for computing with finite simple nonabelian groups (SNAGs)  now entail that periodic Tri-tris is 
capable of {\em finitary universal computation}. This means every function $f:X\rightarrow Y$ for any finite sets $X$ and $Y$ can be realized via an implementation using an encoding into parallel running copies of this Tetris game. 
\begin{theorem}
The periodic 3$\times$4 Tri-tris game is finitarily computationally universal.
\end{theorem}
{\bf Sketch of Proof (construction of~\cite{nehaniv2015symmetry}).}
Let $n = \lceil \log_{60} |X|\rceil $ and $m= \lceil \log_{60} |Y| \rceil$. For each permutation $\pi$ of $Z$ in $A_5$,
fix a particular sequence $w_{\pi}$ of Tetris events 
yielding $\pi$. One encodes distinct members of $X$ each uniquely into $n$ such sequences $(w_{\pi_1},\ldots, w_{\pi_n})$.   Similarly, encode members of $Y$ uniquely in $m$-tuples of permutations $(\pi'_1,\ldots, \pi'_m)$, $\pi'_i \in A_5$.  This yields an encoding of $f$ as a mapping from $n$-tuples of the 60 different $w_{\pi}$ sequences to $m$-tuples of permutations in $A_5$. Now by a theorem of Maurer and Rhodes~\cite{maurer1965property}, each of the $m$ components of the encoded $f$ can be computed by some fixed {\em polynomial expression} over this SNAG. That is, each is some finite concatenation of the fixed sequences $w_{\pi}$ giving permutations in $A_5$ and $n$ free variables which take values in event sequences according to the encoding of $X$ (with repetitions possible).  The evaluation of these polynomial expressions with sequences encoding a member of $X$ substituted in for the $n$ variables consists of running Tri-tris and permutes states in parallel copies of the game (each in a configuration from $Z$). The result in Tri-tris comprises $m$ permutations of $Z$ lying in $A_5$ uniquely encoding the value of $f$ in $Y$.   It suffices to use $5m$ copies of Tri-tris since $|Z|=5$ to determine the $m$ permutations encoding $y=f(x)$ with $x\in X, y\in Y$; actually since we are dealing with permutations $4m$ copies of the game suffice.
\section{Conclusion}

We cast Tetris as a finite transformation semigroup, and show that the complexity of the game grows very quickly with the size of the game board.  Modifying the rules of Tetris to restart on completion yields finite simple nonabelian groups (SNAGs) in the holonomy decomposition.  This entails finite universal computational capacity of periodic variants of Tetris.  While we found computationally that non-periodic Tetris examples had only trivial subgroups in their decompositions, it remains an open problem whether this is the case in all non-periodic variants. It also remains to determine the Krohn-Rhodes complexity and which SNAGs occur in other periodic versions of Tetris.

\section*{Appendix: Krohn-Rhodes Theory and the Holonomy Decomposition}

\begin{figure}
    \centering
    \includegraphics[width = 200px]{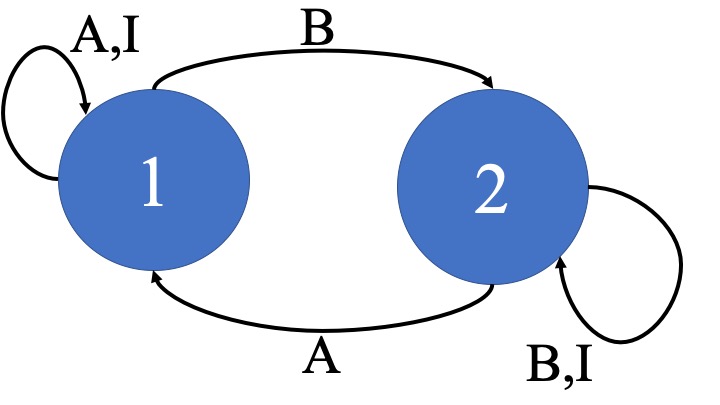}
    \caption{The flip-flop monoid on $ X = \{1,2\}$ is given by the set of transformations $S = \{A,B,I\}$}
    \label{flipflop}
\end{figure}

 The Krohn-Rhodes (KR) theorem describes a general decomposition of transformation semigroups in terms of wreath products of the finite simple groups and the flip-flop monoid. A visualization of the flip-flop monoid is shown in Figure \ref{flipflop}.
 
\begin{theorem}[Krohn-Rhodes decomposition \cite{krohn1965algebraic}]
\label{KR}
A finite transformation semigroup $(X,S)$, with states $X$ and semigroup $S$ acting on states by transformations, has a decomposition 
    $$(X,S) \text{ divides } H_1 \wr H_2 \wr H_3 ... \wr H_n$$
    with components $H_1, H_2, \ldots , H_n$, such that
    each $H_i$ is a finite simple group dividing $S$ or the flip-flop monoid.

\end{theorem}
The decomposition given by this theorem tends to be far from optimal in practice. Therefore, most practical implementations of semigroup decomposition use the holonomy method described in \cite{eilenberg1976automata}. 

We will reproduce the relevant definitions and theorems here. Define the set $Q$ as 
 $$Q = \{\{X \cdot s\} | s \in S\} \cup \{X\} \cup \{\{a\} | a \in X \}$$
 then we can define a relation on $Q$ called subduction.

\begin{definition}[Subduction]
 Let $S^I$ denote $S$ with a new identity element appended.  Given, $A,B \in Q$, we define an reflexive, transitive relation on $Q$, $$A \leq B \iff \exists s \in S^I, A \subseteq B \cdot s$$
Furthermore, let $A < B$ if $A \leq B$ but not $B \leq A$.
 This relation, which we will call \emph{subduction}, induces an equivalence relation on $Q$: \\$A \equiv B \iff A \leq B$ and  $B \leq A$.
      For each equivalence class $A /\!\!\equiv$ in $Q /\!\!\equiv$, let $\bar{A}$ be a unique representative.
\end{definition}
 
 \begin{definition}[Tiles]
Define $A$ to be a \emph{tile} of $B$ if $A \subsetneq B$ and $$\forall Z \in Q,  (A \leq Z \leq B \implies Z = A \mbox{ or } Z = B)$$ 

  If $A \in Q$ with $|A|>1$, the set of tiles of $A$ is $\Theta_A \subset Q$
 \end{definition}

 \begin{definition}[Holonomy group]
 The \emph{holonomy group}, written $H_A$, of $A$ is the set of permutations of $\Theta_A$ induced by the elements of $S^I$.
 If we let $H_A$ act on $\Theta_A$, then $(\Theta_A,H_A)$ is the holonomy permutation group of $A$.
 \end{definition}
  
 \begin{definition}[Height of an Image Set]
  The \emph{height} of $A\in Q$ is $h(A)$, where $h(A)$ is the length of the longest strict subduction chain up to $A$.
 \end{definition}
 
 We are now able to state the holonomy decomposition theorem, which asserts that the semigroup $(X,S)$ divides a cascade product, from which the Krohn-Rhodes (KR) decomposition (Theorem~\ref{KR}) can be derived. The holonomy theorem describes the transformation semigroup $(X,S)$ in terms of symmetries in the way transformations in $S$ act on the set of tiles of the $\bar{A} \in Q$.
 
\begin{theorem}[Holonomy decomposition \cite{eilenberg1976automata}]
\label{holonomy}
 Let $(X,S)$ be a finite transformation semigroup, with $h = h(X)$ the height of $X$. For each $i$ with $1\leq i \leq h$, let $$(\Phi_i, \mathfrak{H_i}) = \prod_{\{A\in Q: h(A)=i,\, \bar{A}=A\}} (\Theta_{\bar{A}}, H_{\bar{A}})$$
   $(\Phi_i, \mathfrak{H_i})$ is a permutation group and $(\Phi_i, \overline{\mathfrak{H_i}})$ is the permutation-reset transformation semigroup obtained by appending all constant maps to $\mathfrak{H_i}$.
Then
$$(X,S) \text{ divides } (\Phi_1, \overline{\mathfrak{H_1}}) \wr (\Phi_2, \overline{\mathfrak{H_2}}) \wr ... \wr (\Phi_h, \overline{\mathfrak{H_h}}).$$

\end{theorem} 
 
\bibliographystyle{spmpsci}
\bibliography{ref}

\end{document}